\documentclass[12pt]{article}
\usepackage{amsmath,amsfonts,amssymb,amsthm,enumerate}

\newcommand{\R}{{\mathbb R}}
\newcommand{\C}{{\mathbb C}}
\newcommand{\Z}{{\mathbb Z}}

\newcommand{\F}{{\mathbb F}}

\newcommand{\Aut}{{\rm Aut}}
\newcommand{\AutO}{{O}}

\newcommand{\eqa}{\begin{eqnarray}}
\newcommand{\eeqa}{\end{eqnarray}}
\newcommand{\eqn}{\begin{eqnarray*}}
\newcommand{\eeqn}{\end{eqnarray*}}

\newcommand{\Hom}{{\rm Hom}}

\newtheorem{dfn}{Definition}[section]
\newtheorem{pro}[dfn]{Proposition}
\newtheorem{thm}[dfn]{Theorem}

\newtheorem{lem}[dfn]{Lemma}
\newtheorem{cor}[dfn]{Corollary}
\newtheorem{rem}[dfn]{Remark}
\newtheorem{note}[dfn]{Note}

\newcommand{\Imm}{{\rm Im\ }}
\newcommand{\Ker}{{\rm Ker\ }}
\newcommand{\NO}{\,{\raise0.25em\hbox{$\mathop{\hphantom{\cdot}}\limits^{_{\circ}}_{^{\circ}}$}}\,}

\newcommand{\qe}{\qed\vskip2ex}
\makeatletter
\@addtoreset{equation}{section}

\makeatother

\def\bl{\begin{lem}\sl}
\def\el{\end{lem}}
\def\bt{\begin{thm}\sl}
\def\et{\end{thm}}
\def\bp{\begin{pro}\sl}
\def\ep{\end{pro}}
\def\br{\begin{rem}\sl}
\def\er{\end{rem}}
\def\bc{\begin{cor}\sl}
\def\ec{\end{cor}}
\def\bd{\begin{dfn}\rm}
\def\ed{\end{dfn}}
\def\bn{\begin{note}\rm}
\def\en{\end{note}}
\def\proof{{\it Proof.}}

\title{\begin{flushright}
\end{flushright}\Large The automorphism groups\\ of the vertex operator algebras $V_L^+$: general case}
\author{Hiroki SHIMAKURA\footnote{The author was supported by the Japan Society for the Promotion of Science Research Fellowships for Young Scientists and COE grant of Hokkaido University.}
}
\date{\small\it Department of Mathematics, Hokkaido University\\
Kita 10, Nishi 8, Kita-Ku, Sapporo, Hokkaido, 060-0810, Japan.\\
{\rm e-mail: shimakura@math.sci.hokudai.ac.jp}
}

\begin{document}
\maketitle

\begin{abstract}
In this article, we give a method of calculating the automorphism groups of the vertex operator algebras $V_L^+$ associated with even lattices $L$.
For example, by using this method we determine the automorphism groups of $V_L^+$ for even lattices of rank one, two and three, and even unimodular lattices.
\end{abstract}


\section*{Introduction}
Let $L$ be a (positive-definite) even lattice and let $V_L^+$ be the fixed-points of the VOA $V_L$ associated with $L$ under an automorphism $\theta_{V_L}$ lifting the $-1$-isometry of $L$.
The automorphism groups $\Aut(V_L^+)$ of the VOAs $V_L^+$ were described in \cite{DG1} for lattices $L$ of rank $1$, in \cite{DG2} for lattices $L$ of rank $2$, and in \cite{Sh2} for lattices $L$ without roots.
The primary purpose of this article is to generalize the method of calculating $\Aut(V_L^+)$ in \cite{Sh2} to all even lattices $L$.

Let $V$ be a VOA and let $G$ be an automorphism group of $V$.
Then the subspace $V^G$ of points fixed by $G$ is a subVOA.
Clearly $N_{\Aut(V)}(G)$ acts on $V^G$.
Then the question arises as to whether or not any automorphism of $V^G$ comes from $N_{\Aut(V)}(G)$.
Take $V$ to be the VOA $V_L$ and $G$ to be the group generated by the involution $\theta_{V_L}$.
Then the quotient group $H_L$ of $C_{\Aut(V_L)}(\theta_{V_L})$ by the subgroup $\langle\theta_{V_L}\rangle$ acts faithfully on $V_L^+$.
In \cite{DG2} it was shown that $\Aut(V_L^+)$ coincides with $H_L$ if $L$ does not have vectors of norm $2$ or $4$ and the rank of $L$ is greater than $1$.
In this article, we can obtain a definitive answer:
$\Aut(V_L^+)$ is larger than $H_L$ if and only if $L$ is obtained by Construction B or is isomorphic to the $E_8$-lattice.

\medskip

We recall the method of \cite{Sh2}.
Let $S_L$ denote the set of all isomorphism classes of irreducible $V_L^+$-modules.
Then $\Aut(V_L^+)$ acts on $S_L$.
It was shown that the stabilizer of the isomorphism class $[0]^-$ of the irreducible $V_L^+$-module $V_L^-$ is equal to the subgroup $H_L$.
The orbit $Q_L$ of $[0]^-$ was determined when $L$ has no roots.
Moreover, $Q_L$ was regarded as a subset of an elementary abelian $2$-group by using the fusion rules of $V_L^+$.
Hence there exists a group homomorphism from $\Aut(V_L^+)$ to a general linear group over $\F_2$.
Then by using the kernel and image $\Aut(V_L^+)$ can be described.

\medskip

The main result of this article is the following:
The orbits $Q_L$ are determined for all even lattices $L$ (Theorem \ref{Torbit}).
This allows us to determine the automorphism group of $V_L^+$.

\medskip

We explain our method of determining $Q_L$.
Since the action of $\Aut(V_L^+)$ on $Q_L$ preserves the graded dimensions and fusion rules, we obtain some necessary conditions satisfied by elements of $Q_L$.
For any element $W$ of untwisted type in $S_L$ satisfying the conditions, we will show that there exists an automorphism exchanging $[0]^-$ and $W$.
To do this, we use a characterization of even lattices obtained by Construction B (Theorem \ref{TM0}) and certain automorphisms given in \cite{FLM}. 
Thus we obtain sufficient and necessary conditions for isomorphism classes of untwisted type to belong $Q_L$.
Moreover we will classify even lattices $L$ such that $Q_L$ contains isomorphism classes of twisted type.
Determining isomorphism classes of twisted type in $Q_L$, we obtain the orbit $Q_L$.

\medskip

Throughout this article, we will work over the field $\C$ of complex numbers unless otherwise stated.
We denote the set of integers by $\Z$ and the rings of integers modulo $p$ by $\Z_p$.
We often identify $\Z_2$ with the field $\F_2$ of two elements.
Let $\Omega_n$ denote the set $\{1,2,\dots,n\}$ for $n\in\Z_{>0}$.
We view the power set of $\Omega_n$ as an $n$-dimensional vector space over $\F_2$ naturally.
For a code $C$ and $l\in\Z$, let $C_l$ denote the set of codewords of $C$ of weight $l$.
For a subset $U$ of an $n$-dimensional vector space $\R^n$ over the real field $\R$ and $m\in\R$, let $U_m$ denote the set of vectors in $U$ of norm $m$.
For a lattice $L$, the dual lattice of $L$ is denoted by $L^*$.
For a group $G$ and its subgroup $H$, $N_G(H)$ and $C_G(H)$ denote the normalizer and centralizer of $H$ in $G$ respectively.
Let $V$ be a VOA and let $(M,Y_M)$ be a $V$-module.
For an automorphism $g$ of $V$, let $M\circ g$ denote the $V$-module $(M,Y_{M\circ g})$ defined by $Y_{M\circ g}(v,z)=Y_M(gv,z)$, $v\in V$.

\medskip

{\it Acknowledgments.} The author would like to thank Professor Atsushi Matsuo for valuable suggestions and helpful advice.
He also thanks Professor Masahiko Miyamoto for useful comments and Professor Toshiyuki Abe for reading the manuscript.

\section{Preliminaries}
In this section, we recall or give some definitions and facts necessary in this article.
\subsection{Construction B}
In this subsection, we recall a standard method for constructing lattices from linear binary codes.

Let $n$ be a positive integer and let $\{\alpha_i|\ i\in\Omega_n\}$ be an orthogonal basis of $\R^n$ satisfying $\langle\alpha_i,\alpha_j\rangle=2\delta_{i,j}$.
For a subset $J\subset\Omega_n$, we set $\alpha_J=\sum_{i\in J}\alpha_i$.
Let $C$ be a binary code of length $n$.
Then 
\begin{eqnarray}
L_B(C)=\sum_{c\in C}\Z\frac{1}{2}\alpha_c+\sum_{i,j\in\Omega_n}\Z(\alpha_i+\alpha_j)\label{Eq:ConstB}
\end{eqnarray}
is called the {\it lattice obtained by Construction B} from $C$.
We note that $L_B(C)$ is even if and only if $C$ is doubly even.
We call $\{\pm\alpha_i|\ i\in\Omega_n\}$ a {\it frame} of $L_B(C)$ with respect to the expression (\ref{Eq:ConstB}).
The following lemma is easy to prove.
\bl\label{LL1} $|{L_{B}(C)}_2|=8|C_4|$.\el

\subsection{Vertex operator algebra $V_L^+$}
In this subsection, we review some properties of the vertex operator algebra $V_L^+$.
For the details of its construction, see \cite{FLM}.

Let $L$ be a (positive-definite) even lattice and let $\hat{L}$ be a central extension:
\begin{eqnarray*}
1\to\langle\kappa_L|\ \kappa_L^2=1\rangle\to\hat{L}\ \bar{\to}\ L\to 1
\end{eqnarray*}
such that $[a,b]=\kappa_L^{\langle\bar{a},\bar{b}\rangle}$ for $a,b\in \hat{L}$.
Let $\theta_{\hat{L}}$ be an involution of $\hat{L}$ induced by the $-1$-isometry of $L$.
Set $K_L=\{a^{-1}\theta_{\hat{L}}(a)|\ a\in \hat{L}\}$.
Then $K_L$ is a normal subgroup of $\hat{L}$.
Let $V_L$ denote the VOA associated with $L$.
The automorphism group $\Aut(V_L)$ of $V_L$ contains an involution $\theta_{V_L}$ induced by $\theta_{\hat{L}}$.
Its fixed-points on $V_L$ is denoted by $V_L^+$.
Then $V_L^+$ is a subVOA of $V_L$.

In \cite{DN2,AD}, it was shown that any irreducible $V_L^+$-module is isomorphic to one of $V_{\lambda+L}^\pm$ $(\lambda\in L^*\cap (L/2))$, $V_{\mu+L}$ $(\mu\in L^*\setminus (L/2))$ and $V_L^{T_\chi,\pm}$, where $T_\chi$ is an irreducible $\hat{L}/K_L$-module with central character $\chi$.
In this article, we use the following notation: $[\mu]$, $[\lambda]^\pm$ and $[\chi]^\pm$ denote the isomorphism classes of $V_{\mu+L}$, $V_{\lambda+L}^\pm$ and $V_{L}^{T_{\chi},\pm}$ respectively.
The isomorphism classes $[\mu]$, $[\lambda]^\pm$ are called {\it untwisted type} and the isomorphism classes $[\chi]^\pm$ are called {\it twisted type}.

\bn In this article, we take an involution on $V_L^{T_\chi}$ induced by the identity operator on $T_\chi$ and consider the $\pm1$-eigenspace $V_L^{T,\pm}$.
However in \cite{FLM} an involution on $V_L^{T_\chi}$ induced by the $-1$-isometry on $T_\chi$ is used.
\en

The fusion rules of $V_L^+$ were determined in \cite{Ab,ADL}.
In particular the following hold.

\bl\label{LFl}{\rm \cite{Ab,ADL}}\sl
\begin{enumerate}
\item Let $\lambda$ be a vector in $L^*\cap (L/2)$.
Then the fusion rules $[0]^-\times[\lambda]^\pm=[\lambda]^\mp$ hold.
\item Let $\lambda$ be a vector in $L^*\cap (L/2)$ satisfying $\langle\lambda,\lambda\rangle\in\Z$.
Then the fusion rule $[\lambda]^\varepsilon\times[\lambda]^\varepsilon=[0]^+$ holds for any $\varepsilon\in\{\pm\}$.
\item Let $W_1$ and $W_2$ be isomorphism classes of irreducible modules of $V_L^+$.
If isomorphism classes of twisted type appear in $W_1\times W_2$ then one of $W_1$ and $W_2$ is of twisted type.
\end{enumerate}
\el

\subsection{Automorphism groups of $V_L$ and $V_L^+$}
In this section, we review the results on automorphism groups of $V_L$ and $V_L^+$ for even lattices $L$.

We start by recalling the automorphism group of $V_L$.
For a lattice $L$, we denote by $O(L)$ the group of automorphisms of $L$ which preserve the bilinear form.
Let $O(\hat{L})$ denote the group of automorphisms of $\hat{L}$ which preserve the bilinear form on the quotient of $\hat{L}$ by its normal subgroup of order $2$.
For $g\in\AutO(\hat{L})$, let $\bar{g}$ denote the linear automorphism of $L$ defined by $\bar{g}(\bar{a})=\overline{g(a)}$, $a\in \hat{L}$.
We view an element $f\in\Hom(L,\Z_2)$ as the automorphism of $\hat{L}$ which sends $a$ to $\kappa_L^{f(\bar{a})}a$.
Hence we obtain an embedding $\Hom(L,\Z_2)\subset \AutO(\hat{L})$.
In Proposition 5.4.1 of \cite{FLM}, the following sequence is exact:
\begin{eqnarray}
1\rightarrow\Hom(L,\Z_2)\hookrightarrow \AutO(\hat{L})\ \bar{\rightarrow}\ \AutO(L)\rightarrow 1.\label{exactauto2}
\end{eqnarray}
In \cite{DN1}, the automorphism group  $\Aut(V_L)$ of $V_L$ was described as follows:
\bp{\rm \cite[Theorem 2.1]{DN1}} Let $L$ be an even lattice.
Then $\Aut(V_L)=N(V_L)O(\hat{L})$, where $N(V_L)=\langle\exp(v_0)|\ v\in (V_L)_1\rangle$ is a normal subgroup of $\Aut(V_L)$.
Moreover, $\Aut(V_L)/N(V_L)$ is isomorphic to a quotient group of $O(L)$.
\ep

In \cite{Do} it was shown that any irreducible $V_L$-module is isomorphic to $V_{\lambda+L}$ for some $\lambda\in L^*$.
The group $\Aut(V_L)$ acts on the set of isomorphism classes of irreducible $V_L$-modules as follows.
\bl\label{RAct} \begin{enumerate}
\item Any element of $N(V_L)$ fixes all isomorphism classes of irreducible $V_L$-module.
\item Let $g$ be an element of $O(\hat{L})$ and let $\lambda$ be a vector in $L^*$.
Then $g$ sends the isomorphism class of $V_{\lambda+L}$ to that of $V_{\bar{g}^{-1}(\lambda)+L}$.
\end{enumerate}
\el

Now, let us consider automorphisms of $V_L^+$.
By the definition of $V_L^+$, the centralizer $C_{\Aut(V_L)}(\theta_{V_L})$ acts on $V_L^+$.
Set $$H_L=C_{\Aut(V_L)}(\theta_{V_L})/\langle\theta_{V_L}\rangle.$$
Then $H_L$ acts faithfully on $V_L^+$, and $H_L\subset\Aut(V_L^+)$.
Let $S_L$ denote the set of all isomorphism classes of irreducible $V_L^+$-modules.
Then $H_L$ is characterized as follows.

\bl\label{LStab} {\rm \cite[Proposition 3.10]{Sh2}} The group $H_L$ is the stabilizer of $[0]^-$ under the action of $\Aut(V_L^+)$ on $S_L$.
\el

Since $\theta_{V_L}$ belongs to the center of $O(\hat{L})$, $H_L$ contains $O(\hat{L})/\langle\theta_{V_L}\rangle$.

\bl\label{LActO}{\rm \cite[Proposition 2.9]{Sh2}} For $g\in\AutO(\hat{L})/\langle\theta_{V_L}\rangle$, we have 
\eqn
[\mu]\circ g&=&[\bar{g}^{-1}[(\mu)],\ \mu\in L^*\setminus (L/2),\\ 
{}\{[\lambda]^{\pm}\circ g\}&=&\{[\bar{g}^{-1}(\lambda)]^\pm\},\ \lambda\in L^*\cap(L/2),\\
{}[0]^{\pm}\circ g&=&[0]^\pm.
\eeqn
Moreover for $\lambda\in L^*\cap(L/2)$ there exists an automorphism $h$ of $V_L^+$ such that $[\lambda]^+\circ h=[\lambda]^-$.
\el

Let $Q_L$ denote the orbit of $[0]^-$ under the action of $\Aut(V_L^+)$ on $S_L$.
Since automorphisms of a VOA preserves the fusion rules and the graded dimensions, we have the following inclusions:

\bl\label{LOrbit}{\rm \cite[Lemma 3.12]{Sh2}} Let $L$ be an even lattice of rank $n$.
\begin{enumerate}
\item If $n\neq8,16$ then $Q_L\subseteq\{[0]^-,[\lambda]^\pm|\ \lambda\in L^*\cap (L/2),\ |(\lambda+L)_2|=2n+|L_2|\}$.
\item If $n=8$ then $Q_L\subseteq\{[0]^-,[\lambda]^\pm,\ [\chi]^-|\ \lambda\in L^*\cap (L/2),\ |(\lambda+L)_2|=2n+|L_2|\}$, where $\chi$ ranges over the central characters of $\hat{L}/K_L$ with $\chi(\kappa K_L)=-1$.
\item If $n=16$ then $Q_L\subseteq\{[0]^-,[\lambda]^\pm,\ [\chi]^+|\ \lambda\in L^*\cap (L/2),\ |(\lambda+L)_2|=2n+|L_2|\}$, where $\chi$ ranges over the central characters of $\hat{L}/K_L$ with $\chi(\kappa K_L)=-1$.
 \end{enumerate}
\el

Recall that a lattice $L$ is said to be {\it $2$-elementary} if $2L^*\subset L$, and said to be {\it totally even} if both $\sqrt2L^*$ and $L$ are even.

\bl\label{L2elem}
\begin{enumerate}
\item{\rm \cite[Proposition 3.14]{Sh2}} If $Q_L$ contains isomorphism classes of twisted type then $L$ is $2$-elementary totally even.
\item{\rm \cite[Lemma 3.6]{Sh2}} If $L$ is $2$-elementary totally even then isomorphism classes of twisted type with the same sign are conjugate under the action of $\Aut(V_L^+)$.
\end{enumerate}
\el

\subsection{Extra automorphisms of $V_L^+$}
In this subsection we review automorphisms of $V_L^+$ not in $H_L$ from \cite{FLM}.

Let $C$ be a doubly even code of length $n$ and let $L$ be the lattice obtained by Construction B from $C$ with frame $\{\alpha_i|\ i\in\Omega_n\}$.
In Chapter 11 of \cite{FLM}, an automorphism $\sigma$ not in $H_L$ was constructed.
This automorphism satisfies $[0]^{-}\circ\sigma=[\alpha_1]^+$.
Lemma \ref{LStab} shows that $\sigma\notin H_L$.
By Lemma \ref{LActO}, there exists an automorphism $h$ of $V_L^+$ such that $[\alpha_1]^+\circ h=[\alpha_1]^-$.

We now assume that $C$ contains the all-ones codeword.
Let us see the action of $\sigma$ on some isomorphism classes of irreducible $V_L^+$-modules.
Set $\beta_n=\alpha_{\Omega_n}/4$ and $\gamma_n=\alpha_{\Omega_n}/4-\alpha_1$.
By the assumption, vectors $\beta_n$ and $\gamma_n$ belong to $L^*\cap (L/2)$.
By \cite[Theorem 11.5.1]{FLM}, $\{[\beta_n]^{\pm}\circ\sigma,[\gamma_n]^{\pm}\circ\sigma \}=\{[\chi_1]^\pm, [\chi_2]^\pm\}$ for some central characters $\chi_i$ of $\hat{L}/K_L$.
Comparing the graded dimensions, we obtain $\{[\beta_n]^{\pm}\circ\sigma\}=\{[\chi_1]^+,[\chi_2]^+\}$ and $\{[\gamma_n]^{\pm}\circ \sigma\}=\{[\chi_1]^-,[\chi_2]^-\}$.

The result is summarized in the following lemmas.

\bl\label{LExt}{\rm \cite{FLM}} Let $L$ be an even lattice obtained by Construction B with frame $\{\alpha_i\}$.
Then the orbit of $[0]^-$ contains $[\alpha_1]^\pm$.
In particular the cardinality of the orbit of $[0]^-$ is grater than $1$.
\el

\bl\label{LExtra}{\rm \cite{FLM}} Let $L$ be the even lattice obtained by Construction B from a doubly even code $C$ containing the all-one codeword.
Let $\varepsilon\in\{\pm\}$.
\begin{enumerate}
\item The orbit of the isomorphism classes $[\beta_n]^\varepsilon$ contains isomorphism classes of twisted type with sign $+$.
\item The orbit of the isomorphism classes $[\gamma_n]^\varepsilon$ contains isomorphism classes of twisted type with sign $-$.
\end{enumerate}
\el

\section{Characterization of even lattices obtained by Construction B}
In this section, we characterize even lattices obtained by Construction B.
We will later use our characterization to determine the automorphism group of $V_L^+$.

Let $L$ be a (positive-definite) even lattice of rank $n$.
We set
\begin{eqnarray}
R_L=\Big\{\lambda+L\in L^*/L\Big|\ \lambda\in L/2,\ |(\lambda+L)_2|\ge2n+|L_2|\Big\}.
\end{eqnarray}

\bn The definition of $R_L$ comes from the necessary conditions satisfied by isomorphism classes of untwisted type in $Q_L$ (cf. Lemma \ref{LOrbit}).
\en

Then even lattices obtained by Construction B are characterized as follows.

\bt\label{TM0} Let $L$ be an even lattice of rank $n$.
Then the following conditions are equivalent:
\begin{enumerate}
\item $L$ is obtained by Construction B.
\item The set $R_L$ is not empty.
\end{enumerate}
\et

To prove this theorem, we need some lemmas.

\bl\label{Llattice1} Let $L$ be an even lattice of rank $n$ and let $\lambda+L$ be an element of $R_L$.
Then $(\lambda+L)_2$ contains an orthogonal basis of $\R^n$.
\el
\proof\ Since $L$ is even and $\lambda\in L^*$, the norms of vectors in $\lambda+L$ are contained in $\langle\lambda,\lambda\rangle+2\Z$.
It follows from $(\lambda+L)_2\neq\phi$ that $\langle\lambda,\lambda\rangle\in2\Z$.
Hence $L^\prime=L+\Z\lambda$ is an even lattice and $L^\prime_2$ forms a root system.
In particular, the inner products of vectors in $(\lambda+L)_2$ are contained in $\{0,\pm1,\pm2\}$.

Let $Y_r=\{y_1,\dots,y_r\}$ be a subset of $(\lambda+L)_2$ such that $\langle y_i,y_j\rangle=2\delta_{i,j}$.
We set $\tilde{Y}_r=\{\pm y|\ y\in Y_r\}$.
Then $|\tilde{Y}_r|=2r$.
To prove this lemma, we will show that if $r<n$ then there exists a vector in $(\lambda+L)_2$ orthogonal to $Y_r$.
Define
\begin{eqnarray*}
X(Y_r)=\Bigl\{x\in(\lambda+L)_2\Big\arrowvert\ \langle x,y_i\rangle\in\{\pm1\}\ {\rm for}\ {}^\exists i\in\Omega_r\Bigr\},
\end{eqnarray*}
where $\Omega_r=\{1,2,\dots,r\}$.
Clearly $\tilde{Y}_r\cap X(Y_r)=\phi$.
For $x\in X(Y_r)$, we set 
\begin{eqnarray}
m(x)=\min\Bigl\{i\in\Omega_r\Big\arrowvert\ \langle y_i,x\rangle\in\{\pm1\}\Bigr\}.\label{Def:m(x)}
\end{eqnarray}
Set $X(Y_r)^+=\{x\in X(Y_r)|\ \langle x,y_{m(x)}\rangle=1\}$.
Then $|X(Y_r)|=2|X(Y_r)^+|$.
Since $|x-y_{m(x)}|^2=|x|^2+|y_{m(x)}|^2-2\langle x,y_{m(x)}\rangle=2$ for $x\in X(Y_r)^+$, we consider the map $\rho :X(Y_r)^+\to \{\{\pm v\}|\ v\in L_2\}$, $x\mapsto\{\pm(x-y_{m(x)})\}$.

Let us show that $\rho$ is injective.
First we suppose $x-y_{m(x)}=x^\prime-y_{m({x^\prime})}$.
If $m(x)=m({x^\prime})$ then $x=x^\prime$.
So we may assume that $m(x)<m(x^\prime)$.
By the definition of $Y_r$ and (\ref{Def:m(x)})
\begin{eqnarray}
\langle y_{m(x)},y_{m(x^\prime)}\rangle=\langle x^\prime,y_{m(x)}\rangle=0.\label{Eq:La1}
\end{eqnarray}
Hence we have 
\begin{eqnarray*}
2=\langle x-y_{m(x)},x^\prime-y_{m({x^\prime})}\rangle=\langle x,x^\prime-y_{m({x^\prime})}\rangle.
\end{eqnarray*}
Since both $x$ and $x^\prime-y_{m(x^\prime)}$ belong to $L^\prime_2$, we have $x=x^\prime-y_{m({x^\prime})}$.
However it contradicts $(\lambda+ L)\cap L=\phi$ since $x\in \lambda+L$ and $x^\prime-y_{m(x^\prime)}\in L$.

Next we suppose $x-y_{m(x)}=y_{m({x^\prime})}-x^\prime$.
If $m(x)=m({x^\prime})$ then $x+x^\prime=2(y_{m(x)})$ and $|x|^2=|x^\prime|^2=|y_{m(x)}|^2=2$.
Hence $x=x^\prime=y_{m(x)}$, which is a contradiction.
So we may assume $m(x)< m({x^\prime})$.
Then by (\ref{Eq:La1}), we have
\begin{eqnarray*}
2=\langle x-y_{m(x)},y_{m({x^\prime})}-x^\prime\rangle=\langle x,y_{m(x^\prime)}-x^\prime\rangle,
\end{eqnarray*}
which implies that $x=y_{m(x^\prime)}-x^\prime$.
However, it contradicts $(\lambda+L)\cap L=\phi$.
Hence $\rho$ is injective.
This shows that $|X(Y_r)^+|\le|L_2|/2$, namely $|X(Y_r)|\le|L_2|$.
Since $\tilde{Y}_r\cap X(Y_r)=\phi$, we have $|\tilde{Y}_r\cup X(Y_r)|\le2r+|L_2|$.
So, if $r< n$ then there exists $x\in (\lambda+L)_2$ such that $x\notin X(Y_r)\cup \tilde{Y}_r$, namely $\langle x,Y_r\rangle=0$.
Therefore $(\lambda+L)_2$ contains an orthogonal basis of $\R^n$.
\qe

\bl\label{PM0} Let $L$ be an even lattice of rank $n$.
For $\lambda+L\in R_L$, there exist a doubly even code $C$ and a frame in $\lambda+L$ such that $L=L_B(C)$.
\el

\proof\ By Lemma \ref{Llattice1}, $\lambda+L$ contains vectors $e_1,e_2,\dots, e_n$ satisfying $\langle e_i,e_j\rangle=2\delta_{i,j}$.
Since $2\lambda\in L$, we have $e_i\pm e_j\in L$.
Set $E=\oplus_{i=1}^n \Z e_i$ and $L^\prime=L+\Z \lambda$.
Then $L^\prime/E$ is a subspace of $E^*/E\cong \Z_2^n$.
So we regard $L^\prime/E$ as a binary code $C$ of length $n$.
We can choose a basis $B$ in $\{\pm e_i|\ i\in\Omega_n\}$ so that $L$ is the lattice obtained by Construction B from $C$ with frame $B$ (cf. the proof of \cite[Proposition 1.8]{Sh2}).
\qe 

{\it Proof of Theorem \ref{TM0}}.
Suppose (1).
Let $\{\alpha_i\}$ be the frame.
Then $|(\alpha_1+L)_2|=2n+|L_2|$, and $\alpha_1+L\in R_L$.
Hence $(1)\Rightarrow (2)$.
It follows from Lemma \ref{PM0} that $(2)\Rightarrow(1)$.\qe

\br The proof of Theorem \ref{TM0} implies that $|(\lambda+L)_2|=2n+|L_2|$ for any $\lambda+L\in R_L$.
\er

Let us show some lemmas by using Theorem \ref{TM0}.
Let $L$ be the even lattice obtained by Construction B from a doubly even code $C$ of length $n$ with frame $\{\pm\alpha_i|\ i\in\Omega_n\}$.
Set $\beta_n=\alpha_{\Omega_n}/4$ and $\gamma_n=\alpha_{\Omega_n}/4-\alpha_1$.

\bl\label{Lgamma} The following conditions are equivalent:
\begin{enumerate}
\item $\gamma_n+L\in R_L$.
\item $n=8$ and $C$ contains the all-one codeword.
\end{enumerate}
\el
\proof\ Suppose (2).
Since $C$ contains the all-one codeword, $\gamma_8\in L/2$.
Clearly $\gamma_8\in L^*$.
It is easy to see that $$(\gamma_8+L)_2=\{\pm(\alpha_{\Omega_8}/4-\alpha_i), \alpha_{\Omega_8}-\alpha_c/2+\alpha_j,\alpha_{\Omega_8}-\alpha_c/2-\alpha_k|\  c\in C_4, i\in\Omega_8,j\in c,k\in\Omega_8\setminus c\}.$$
Hence $|(\gamma_8+L)_2|=16+8|C_4|=16+|L_2|$ by Lemma \ref{LL1}.
Thus $\gamma_8+L\in R_L$.

Conversely, we suppose (1).
Since the norm of $\gamma_n$ is minimal in $\gamma_n+L$ and it is $1+n/8$, the rank $n$ of $L$ must be $8$.
Since $\gamma_8\in L/2$, we obtain $\alpha_{\Omega_8}/2\in L$.
Hence the all-one codeword belongs to $C$.
\qed

\bl\label{Lbeta} The following conditions are equivalent:
\begin{enumerate}
\item $\beta_n+L\in R_L$.
\item $n=16$ and $C$ contains a subcode isomorphic to the Reed-Muller code $RM(1,4)$.
\end{enumerate}
\el
\proof\ Let $k$ be the dimension of $C$.
Suppose (2).
In \cite{PLF}, doubly even codes of length $16$ containing the all-one codeword were classified.
In particular doubly even codes of length $16$ containing $RM(1,4)$ can be classified.
Hence we obtain $|C_4|=0,4,12,28$ for $k=5,6,7,8$ respectively.
So $32+|L_2|=2^k$.
On the other hand, 
\begin{eqnarray}
(\beta_{16}+L)_2=\{\beta_{16}-\alpha_c/2|\ c\in C\}.\label{Eq:beta2}
\end{eqnarray}
Hence $|(\beta_{16}+L)_2|=2^k$.
Therefore $\beta_{16}+L\in R_L$.

Conversely we suppose (1).
Since the norm of $\beta_{n}$ is minimal in $\beta_{n}+L$ and it is $n/8$, the rank $n$ of $L$ must be $16$.
By Lemma \ref{Llattice1}, $(\beta_{16}+L)_2$ contains an orthogonal basis $F$.
Set $\tilde{F}=\{\pm v|\ v\in F\}$.
By (\ref{Eq:beta2}), $\tilde{F}=\{\beta_{16}-\alpha_c/2|\ c\in D\}$ for some subset $D$ of $C$.
Clearly $|D|=32$.
Let $d$ be an element of $D$.
Set $D^0=d+D$.
Then $\tilde{F}^0=\{\beta_{16}-\alpha_c/2|\ c\in D^0\}$ is a set of $32$ vectors, two of which are equal, opposite, or orthogonal.
Since $D^0$ contains the all-zero codeword, $D^0$ consists of the all-zero and all-one codewords and $30$ codewords with weight $8$.
Moreover, the cardinality of any intersection of codewords with weight $8$ in $D^0$ is $0$, $4$ or $8$.
Hence $D^0$ must be isomorphic to the Reed-Muller code $RM(1,4)$.
Therefore $C$ contains a subcode isomorphic to the Reed-Muller code $RM(1,4)$.
\qed

\section{Automorphism groups of $V_L^+$ for even unimodular lattices of rank $8$ and $16$}\label{SAut816}
In this section, we determine the automorphism groups of $V_L^+$ for even unimodular lattices of rank $8$ and $16$.
In particular, we will compare $\Aut(V_L^+)$ with its subgroup $H_L\cong C_{\Aut(V_L)}(\theta_{V_L})/\langle\theta_{V_L}\rangle$.

Let $L$ be an even unimodular lattice of rank $8$ or $16$.
Then the VOA $V_L^+$ has exactly $4$ isomorphism classes of irreducible $V_L^+$-modules $[0]^\pm$ and $[\chi_0]^\pm$, where $\chi_0$ is the unique faithful character of $\hat{L}/K_L$.
Since the graded dimensions of $[\chi_0]^+$ and $[\chi_0]^-$ are different, the cardinality of the orbit $Q_L$ of $[0]^-$ is $1$ or $2$.
In the following subsections, we will determine $|Q_L|$.

\subsection{Automorphism group of $V_L^+$ for the even unimodular lattice of rank $8$}\label{SAut8}
In this subsection, we study the automorphism group of $V_{E_8}^+$, where $E_8$ is the unique even unimodular lattice of rank $8$ up to isomorphism.

\bl\label{LE8} There are automorphisms of $V_{E_8}^+$ mapping $[0]^-$ to $[\chi_0]^-$.
In particular $Q_{E_8}$ contains isomorphism classes of twisted type.
\el

\proof\ The degree $1$ subspace of $V_{E_8}^+$ forms the simple Lie algebra of type $D_8$ under the $0$-th product and $V_{E_8}^+\cong V_{D_8}$.
By \cite{Do}, $V_{D_8}$ has exactly $4$ non-isomorphic irreducible modules $V_{\lambda+D_8}$, $\lambda+D_8\in D_8^*/D_8$.

On the other hand, there exists an involution $\tau$ of the root lattice of type $D_8$ such that $\tau$ exchanges two elements of $D_8^*/D_8$.
By Lemma \ref{RAct} (2), lifts of $\tau$ exchange two isomorphism classes of irreducible $V_{D_8}$-modules.
This shows that there are automorphisms of $V_{E_8}^+$ mapping $[0]^-$ to $[\chi_0]^-$.
\qed

\bp\label{PE8} The group $H_{E_8}$ is a normal subgroup of $\Aut(V_{E_8}^+)$ of index $2$.
In particular $\Aut(V_{E_8}^+)/H_{E_8}\cong \Z_2$.
\ep
\proof\ Lemma \ref{LE8} shows that the cardinality of the orbit $Q_{E_8}$ of $[0]^-$ is $2$.
By Lemma \ref{LStab} $H_{E_8}$ is a subgroup of the index $2$ of $\Aut(V_{E_8}^+)$.\qe

\subsection{Automorphism groups of $V_L^+$ for even unimodular lattices of rank $16$}
In this subsection, we study the automorphism groups of $V_{L}^+$ for even unimodular lattices $L$ of rank $16$.
It is known that $E_8{\oplus}E_8$ and $\Gamma_{16}$ are the only even unimodular lattices of rank $16$ up to isomorphism (cf. \cite{CS}).
We note that the root sublattice of $\Gamma_{16}$ is of type $D_{16}$.

Let $U$ be a root lattice of type $D_8\oplus D_8$.
Let $N$ be an even overlattice of $U$ such that $|N:U|=2$ and $N_2=U_2$.
It is easy to check that $N$ is unique up to isomorphism.
Since the determinant of $N$ is $4$, there are three unimodular overlattices of $N$.
In particular even unimodular lattices $E_8{\oplus}E_8$ and $\Gamma_{16}$ are obtained as overlattices of $N$.

\bl\label{LE82} Any element of $\Aut(V_{N})$ of $V_N$ fixes all isomorphism class of irreducible $V_N$-modules.
\el
\proof\ The automorphism group $O(N)$ of $N$ acts on $N^*/N$.
Since unimodular overlattices of $N$ are non-isomorphic, $O(N)$ fixes all elements of $N^*/N$.
By Lemma \ref{RAct} (2), we obtain this lemma.
\qe

\bl\label{LD16} The VOAs $V_{\Gamma_{16}}^+$ and $V_{E_8\oplus E_8}^+$ are isomorphic to $V_N$.
\el

\proof\ First we consider the lattice $E_8\oplus E_8$.
Since $V_{E_8}^+$ is isomorphic to $V_{D_8}$, $V_{E_8\oplus E_8}^+$ contains a subVOA $V$ isomorphic to $V_{U}$.
Since $V_{U}$ is rational, $V_{E_8\oplus E_8}^+= V\oplus M$ as $V$-module for some $V$-module $M$.
By the classification of irreducible modules of $V_U$ (\cite{Do}), $M$ is isomorphic to the irreducible $V_U$-module $V_{\lambda+U}$, where $\lambda+U\in U^*/U$ satisfying $\langle\lambda,U\rangle= N$.
Since $V_N=V_{U}\oplus V_{\lambda+U}$ is a simple current extension of $V_U$, $V_{E_8\oplus E_8}^+$ has a unique VOA structure extending its $V$-module structure (cf. Proposition 5.3 in \cite{DM}) and $V_{E_8\oplus E_8}^+\cong V_N$.

Next, we consider the lattice $\Gamma_{16}$.
Since the root sublattice of $\Gamma_{16}$ is $D_{16}$, $V_{\Gamma_{16}}^+$ contains $V_{D_{16}}^+$.
The degree $1$ subspace of $V_{D_{16}}^+$ forms a simple Lie algebra of type $D_{8}\oplus D_8$ under the $0$-th product and $V_{D_{16}}^+\cong V_U$.
Similarly to the argument above, we obtain $V_{\Gamma_{16}}^+\cong V_N$.\qe

This lemma shows that the cardinality of the orbit $Q_{L}$ of $[0]^-$ is $1$ for any even unimodular lattice $L$ of rank $16$.
By Lemma \ref{LStab} $\Aut(V_{L}^+)$ is coincides with $H_L$.

\bp\label{C16} The automorphism group $\Aut(V_{L}^+)$ of $V_L^+$ coincides with $H_L$ for any even unimodular lattice $L$ of rank $16$.
\ep

\section{The orbit of the isomorphism class of $V_L^-$}
In this section, we determine the orbit $Q_L$ of $[0]^-$.
We note that $Q_L$ was determined in \cite{Sh2} when $L$ has no roots.

\bl\label{LMT1} The orbit $Q_L$ contains the isomorphism class $[\lambda]^\varepsilon$ for any $\lambda\in R_L$, $\varepsilon\in\{\pm\}$.
\el
\proof\ By Lemma \ref{LOrbit} any isomorphism class of untwisted type in $Q_L$ must be $[\lambda]^\varepsilon$ for some $\lambda\in R_L$ and $\varepsilon\in\{\pm\}$.
Conversely by Lemma \ref{LExt} and \ref{PM0} $Q_L$ contains $[\lambda]^\varepsilon$ for all $\lambda+L\in R_L$ and $\varepsilon\in\{\pm\}$.\qe

So let us discuss the cases where $Q_L$ contains isomorphism classes of twisted type.
We consider the following three conditions on even lattices $L$:
\begin{enumerate}[(a)]
\item $L$ is obtained by Construction B from a doubly even code of length $8$ containing the all-one codeword.
\item $L$ is obtained by Construction B from a doubly even code of length $16$ containing a subcode isomorphic to the first order Reed-Muller code $RM(1,4)$ of length $16$.
\item $L$ is isomorphic to the $E_8$-lattice.
\end{enumerate}

\bp\label{PTWIST} The orbit $Q_L$ contains isomorphism classes of twisted type if and only if $L$ satisfies (a), (b) or (c).
\ep
\proof\ By Lemma \ref{LExtra}, \ref{Lgamma}, \ref{Lbeta} and \ref{LE8} if $L$ satisfies (a), (b) or (c) then $Q_L$ contains isomorphism classes of irreducible $V_L^+$-modules of twisted type.

So we suppose that $Q_L$ contains an isomorphism class $[\chi]^\varepsilon$ of twisted type.
Then by Lemma \ref{LOrbit} the rank of $L$ is $8$ or $16$, and $\varepsilon=-$ and $+$ if $n=8$ and $16$ respectively.
Moreover by Lemma \ref{L2elem} (1) $L$ is $2$-elementary totally even.
If $L$ is unimodular then $L$ is isomorphic to one of $E_8$, $E_8\oplus E_8$ and $\Gamma_{16}$.
By the result of the previous section, $L$ must be isomorphic to $E_8$.
Hence $L$ satisfies (c).

We now assume that $L$ is not unimodular.
Let us show that $Q_L$ contains $[\lambda]^\delta$ for some $\lambda\in L^*\cap(L/2)$, $\delta\in\{\pm\}$.
By comparing the coefficients of $q$ in the graded dimensions of $V_L^-$ and $V_L^{T_\chi,\varepsilon}$, the theta series of $L$ is written by the Dedekind-eta series.
By using the transformation formula on theta series of lattices and their dual lattices, we can describe the theta series of $L^*$.
In particular, $L^*\setminus L$ has vectors of norm $2$ (cf. the proof of Proposition 3.14 in \cite{Sh2}).
Let $\lambda$ be an element of $L^*$ such that $(\lambda+L)_2\neq\phi$.
Let $g$ be an element of $\Aut(V_L^+)$ such that $[0]^{-}\circ g=[\chi]^\varepsilon$.
By Lemma \ref{LFl} (1), we obtain $[\chi]^\varepsilon\times ([\lambda]^{+}\circ g)=[\lambda]^{-}\circ g$.
By Lemma \ref{LFl} (3) one of $[\lambda]^{\pm}\circ g$ must be of twisted type.
By comparing the graded dimensions, it has the same sign $\varepsilon$.
By Lemma \ref{L2elem} (2), $Q_L$ contains $[\lambda]^\delta$ for some $\delta\in\{\pm\}$.
By Lemma \ref{LOrbit} $\lambda+L\in R_L$.
Thus $L$ is obtained by Construction B from a code $C$ by Theorem \ref{TM0}.

Since $L$ is $2$-elementary totally even, $C$ contains the all-one codeword.
Hence (a) holds if the rank of $L$ is $8$.
Consider the case where $n=16$.
Since the theta series of $L$ is described in terms of the weight enumerator of $C$, we can describe the weight enumerator of $C$.
By using the classification of even codes of length $16$ \cite{PLF}, (b) holds if the rank of $L$ is $16$.\qe

By Lemma \ref{LOrbit}, \ref{L2elem}, \ref{LMT1} and Proposition \ref{PTWIST}, the orbit $Q_L$ is determined.

\bt\label{Torbit} Let $L$ be an even lattice of rank $n$.
\begin{enumerate}
\item If $L$ satisfies (a) or (c) then $Q_L=\{[0]^-,[\lambda]^\pm,\ [\chi]^-|\ \lambda\in L^*\cap (L/2),\ |(\lambda+L)_2|=2n+|L_2|\}$, where $\chi$ ranges over the central characters of $\hat{L}/K_L$ with $\chi(\kappa K_L)=-1$ and $\varepsilon=+$.
\item If $L$ satisfies (b) then $Q_L=\{[0]^-,[\lambda]^\pm,\ [\chi]^+|\ \lambda\in L^*\cap (L/2),\ |(\lambda+L)_2|=2n+|L_2|\}$, where $\chi$ ranges over the central characters of $\hat{L}/K_L$ with $\chi(\kappa K_L)=-1$.
\item If $L$ does not satisfy neither (a), (b) nor (c) then $Q_L=\{[0]^-,[\lambda]^\pm|\ \lambda\in L^*\cap (L/2),\ |(\lambda+L)_2|=2n+|L_2|\}$.
\end{enumerate}
\et

By Lemma \ref{LStab}, Theorem \ref{TM0} and \ref{Torbit}, we have the following corollary.
\bc\label{MT1}\sl The automorphism group $\Aut(V_L^+)$ of $V_L^+$ is greater than $H_L$ if and only if the even lattice $L$ satisfies one of the following:
\begin{enumerate}[{\rm (1)}]
\item $L$ is obtained by Construction B.
\item $L$ is isomorphic to the $E_8$-lattice.
\end{enumerate}
\ec

\section{A method of determining of the shape of the automorphism group of $V_L^+$}
In this section, we give a method of determining the shape of $\Aut(V_L^+)$ for an arbitrary even lattice $L$.
This method is a generalization of that in \cite[Section 3.4]{Sh2}.
For the conditions (a), (b) and (c), see the previous section.
If $L$ satisfies (c) then $\Aut(V_L^+)$ is determined in Section \ref{SAut8}.

First we consider the lattice $L$ satisfying neither (a), (b) nor (c).
Then $Q_L=\{[0]^-, [\lambda]^\pm|\ \lambda\in R_L\}$, and $P_L=\{[0]^+\}\cup Q_L$ has an elementary abelian $2$-group structure under the fusion rules (cf. \cite[Proposition 3.17]{Sh2}).
So we obtain a group homomorphism $\varphi_L$ from $\Aut(V_L^+)$ to $GL(P_L)$.
Since the kernel of $\varphi_L$ is a subgroup of $H_L$, it can be determined.
Moreover the index of $\varphi(H_L)$ in $\Imm\varphi_L$ is equal to the cardinality of $Q_L$.
Hence we can determine the image of $\varphi_L$.
Therefore we can calculate the shape of $\Aut(V_L^+)$ in principle.

Suppose that $L$ satisfies (a) or (b).
In this case, we consider the set $S_L$ of all isomorphism classes of irreducible $V_L^+$-modules.
Since $L$ is $2$-elementary totally even, $S_L$ has an elementary abelian $2$-group structure under the fusion rules (cf. \cite{Ab,ADL} and \cite[Proposition 3.4]{Sh2}).
Moreover $S_L$ has a natural quadratic form associated with a non-singular symplectic form preserved by the action of $\Aut(V_L^+)$ (cf. \cite[Theorem 3.8]{Sh2}).
Hence we obtain a group homomorphism $\psi_L$ from $\Aut(V_L^+)$ to the orthogonal group $O(S_L)$ associated with the quadratic form.
Similarly to the case above, we can determine the image and kernel of $\psi_L$, and we can describe the shape of $\Aut(V_L^+)$ in principle.

\bn For many important lattices $L$ without roots, the shapes of $\Aut(V_L^+)$ were determined in \cite[Section 4]{Sh2} by using this method.
\en

\section{Automorphism groups of VOSAs $V_L^+$ for odd lattices}
Let $L$ be an odd lattice.
In this section, we consider the vertex operator superalgebra $V_L^+$. 
For $i\in\{0,1\}$, set $L^i=\{v\in L|\ \langle v,v\rangle\equiv i\ (2)\}$.
Then $L^0$ is an even sublattice of $L$.
We will describe $\Aut(V_{L}^+)$ by using $\Aut(V_{L^0}^+)$.

Let $\Aut(V_{L^0}^+;V_{L^1}^+)$ denote the subgroup of $\Aut(V_{L_0}^+)$ fixing the isomorphism class of $V_{L^1}^+$.
Let $\alpha$ be a vector in $L^1$.
Then $2\alpha\in L^0$ and $\langle\alpha,\alpha\rangle\in\Z$.
By Lemma \ref{LFl} (2) $[\alpha]^+\times[\alpha]^+=[0]^+$.
Let $\tau$ denote the involution acting as $(-1)^i$ on $V_{L^i}^+$.
Applying \cite[Theorem 3.3]{Sh2} to our case, we obtain $C_{\Aut(V_L^+)}(\tau)/\langle\tau\rangle\cong \Aut(V_{L^0}^+;V_{L^1}^+)$.

On the other hand, any automorphism of $V_L^+$ preserves both $V_{L^0}^+$ and $V_{L^1}^+$ since the graded dimensions of $V_{L^0}^+$ and $V_{L^1}^+$ are in $\Z[[q]]$ and in $\Z q^{1/2}[[q]]$ respectively.
Hence $C_{\Aut(V_L^+)}(\tau)=\Aut(V_L^+)$.
Therefore we have the following proposition.

\bp\label{TVOSA}\sl Let $L$ be an odd lattice.
Then $\Aut(V_L^+)\cong \langle\tau\rangle.\Aut(V_{L^0}^+;V_{L^1}^+)$.
\ep

Since the shape of $\Aut(V_{L^0}^+)$ can be described by the method given in the previous section, $\Aut(V_L^+)$ can be determined in principle.

\section{Examples}
In this section, we calculate $\Aut(V_L^+)$ for some lattices by using the method of Section 5.

\subsection{Even lattices of rank one, two and three}
In this section, we determine $\Aut(V_L^+)$ for even lattices of rank one two and three.

Let $L$ be an even lattice $L$ of rank $n$.
Suppose that $n\le3$.
By Theorem \ref{Torbit}, $Q_L=\{[0]^-,[\lambda]^\pm|\ \lambda\in R_L\}$.
So we consider $R_L$.
By Theorem \ref{TM0}, let us consider even lattices obtained by Construction B.
It is easy to see that a code $C$ of length $n$ is doubly even if and only if $C$ consists of the all-zero codeword.
Hence $L$ is obtained by Construction B if and only if $L\cong 2A_1, \sqrt2(A_1\oplus A_1)$ or $\sqrt2A_3$.
If $L$ is not obtained by Construction B then $\Aut(V_L^+)\cong C_{\Aut(V_L)}(\theta_{V_L})/\langle\theta_{V_L}\rangle$.
The case where $L\cong\sqrt2A_3$ was done in Theorem 4.3 of \cite{Sh2}.
So let us consider the automorphism groups of $V_L^+$ for $2A_1$ and $\sqrt2(A_1\oplus A_1)$.

First we consider the case where $L\cong 2A_1$.
Let $\gamma$ be a generator of $L$.
Then $R_L=\{\gamma/2+L\}$.
Hence $Q_L=\{[0]^-,[\gamma/2]^\pm\}$.
Set $P_L=\{[0]^+\}\cup Q_L$.
Then $P_L$ has an elementary abelian $2$-group structure under the fusion rules and, $P_L\cong \F_2^2$.
So we obtain a group homomorphism $\varphi_L:\Aut(V_L^+)\to GL(P_L)$.
On the other hand, $H_L\cong\Z_2$ and its generator exchanges $[\gamma/2]^+$ and $[\gamma/2]^-$.
Since $\Ker \varphi_L$ is a subgroup of $H_L$, $\varphi_L$ is injective.
Clearly $\varphi_L(H_L)$ is a maximal subgroup of $GL(P_L)\cong S_3$.
Since $\Aut(V_L^+)$ contains automorphisms not in $H_L$ (cf. Lemma \ref{LExt}), $\varphi_L$ is surjective.
Thus we obtain $\Aut(V_L^+)\cong S_3$.

Next let us consider the case where $L\cong \sqrt2(A_1\oplus A_1)$.
Let $\{a_1,a_2\}$ be a basis of $L$ satisfying $\langle a_i,a_j\rangle=4\delta_{i,j}$.
Set $a_i^*=a_i/4$ and $b=2(a_1^*+a_2^*)$.
Then $\{a_1^*,a_2^*\}$ is a basis of the dual lattice of $L$ and $R_L=\{b+L\}$.
So $Q_L=\{[0]^-,[b]^\pm\}$.
Set $P_L=\{[0]^+\}\cup Q_L$.
Then $P_L$ has an elementary abelian $2$-group structure under the fusion rules and $P_L\cong \F_2^2$.
So we obtain a group homomorphism $\varphi_L:\Aut(V_L^+)\to GL(P_L)$.
On the other hand, $H_L$ is isomorphic to the direct product of the dihedral group of order $8$ and the group of order $2$.
The kernel of $\varphi_L$ is isomorphic to $2^3$, and $H_L$ contains elements exchanging $[b]^+$ and $[b]^-$.
So $\varphi_L(H_L)$ is a maximal subgroup of $GL(P_L)\cong S_3$.
Since $\Aut(V_L^+)$ contains automorphisms not in $H_L$, $\varphi_L$ is surjective.
Therefore we obtain $\Aut(V_L^+)\cong 2^3.S_3$.
It is easy to check that $\Aut(V_L^+)\cong S_4\times \Z_2$.

The result is summarized in the following proposition.

\bp\label{TAVL+R1} Let $L$ be an even lattice of rank one, two or three.
Then 
\begin{eqnarray*}
\Aut(V_L^+)\cong\left\{\begin{array}{cc}
 \mbox{$S_3$} & \mbox{${\rm if}\ L\cong 2A_1$},\\
 \mbox{$S_4\times \Z_2$} & \mbox{${\rm if}\ L\cong \sqrt2(A_1\oplus A_1)$},\\
 \mbox{$(2^2:S_4).S_3$} & \mbox{${\rm if}\ L\cong\sqrt2A_3$},\\
 \mbox{$C_{\Aut(V_L)}(\theta_{V_L})/\langle\theta_{V_L}\rangle$} & \mbox{{\rm otherwise}}.
 
 \end{array}
\right.
\end{eqnarray*}
\ep

\bn\label{RAVL+R1}The automorphism groups of $V_L^+$ for lattices of rank one and two were already obtained in \cite{DG1} and \cite{DG2} respectively by using the action of $\Aut(V_L^+)$ on certain homogeneous subspaces of $V_L^+$.
In the articles, more precise structures of $\Aut(V_L^+)$ were described.
\en

\subsection{Even unimodular lattices}
Let $L$ be an even unimodular lattice.
Since the determinant of any lattice obtained by Construction B is not $1$, $L$ is not obtained by Construction B.
Hence $R_L=\phi$ by Theorem \ref{TM0}.
By Theorem \ref{Torbit}, $|Q_L|=2$ if $L\cong E_8$, and $|Q_L|=1$ if $L\not\cong E_8$.
By Lemma \ref{LStab} and Proposition \ref{PE8}, we obtain the following proposition.

\bp Let $L$ be an even unimodular lattice of rank $n$.
Then
\begin{eqnarray*}
\Aut(V_L^+)\cong\left\{\begin{array}{cc}
 \mbox{$(C_{\Aut(V_L)}(\theta_{V_L})/\langle\theta_{V_L}\rangle).\Z_2$} & \mbox{${\rm if}\ {\rm rank} L=8$},\\
 \mbox{$C_{\Aut(V_L)}(\theta_{V_L})/\langle\theta_{V_L}\rangle$} & \mbox{${\rm if}\ {\rm rank} L\ge16$}.
 \end{array}
\right.
\end{eqnarray*}
\ep


\begin{thebibliography}{100000}

\bibitem[Ab]{Ab}
T.\ Abe, Fusion rules for the charge conjugation orbifold, {\it J. Algebra}, {\bf 242} (2001), 624--655.

\bibitem[AD]{AD}
T.\ Abe and C.\ Dong, Classification of irreducible modules for the vertex operator algebra $V\sp +\sb L$: general case. {\it J. Algebra} {\bf 273} (2004), 657--685

\bibitem[ADL]{ADL}
T.\ Abe, C.\ Dong and H.\ Li, Fusion rules for the vertex operator algebras $M(1)^+$ and $V_L^+$, {\it Comm. Math. Phys.} {\bf 253} (2005), 171--219.


\bibitem[CS]{CS}
J.H.\ Conway and N.J.A. Sloane, Sphere packings, lattices and groups, 3rd Edition, Springer, New York, 1999.

\bibitem[Do]{Do}
C-Y.\ Dong,  Vertex algebras associated with even lattices, {\it J.\ Algebra} {\bf 160} (1993), 245--265.

\bibitem[DG1]{DG1}
C.\ Dong and R.L.\ Griess, Rank one lattice type vertex operator algebras and their automorphism groups, {\it J.\ Algebra} {\bf 208} (1998), 262--275.

\bibitem[DG2]{DG2}
C.\ Dong and R.L.\ Griess, The rank two lattice type vertex operator algebras $V_L^+$ and their automorphism groups, math.QA/0409409, preprint.

\bibitem[DM]{DM}
C.\ Dong and G.\ Mason, Rational vertex operator algebras and the effective central charge, {\it Int. Math. Res. Not.}  {\bf 56} (2004), 2989--3008.

\bibitem[DN1]{DN1}
C-Y.\ Dong and K.\ Nagatomo, Automorphism groups and Twisted modules for lattice Vertex operator algebras, {\it Comtemp.\ Math.} {\bf 248} (1999), 117--133

\bibitem[DN2]{DN2}
C-Y.\ Dong and K.\ Nagatomo, Representations of vertex operator algebra $V_L^+$ for rank one lattice $L$, {\it Comm.\ Math.\ Phys.} {\bf 202} (1999), 169--195.

\bibitem[FHL]{FHL}
I.\ Frenkel, Y.\ Huang, J.\ Lepowsky, On axiomatic approaches to vertex operator algebras and modules,  Mem. Amer. Math. Soc. {\bf 104} (1993).

\bibitem[FLM]{FLM}
I.\ Frenkel, J.\ Lepowsky and A.\ Meurman, Vertex operator algebras and the Monster, Pure and Appl.\ Math., Vol.134, Academic Press, Boston, 1989.


\bibitem[PLF]{PLF}
V. Pless, J.S. Leon, and J. Fields, All $Z\sb 4$ codes of type II and length $16$ are known, {\it J.\ Combin.\ Theory Ser.\ A} {\bf 78} (1997), 32--50.


\bibitem[Sh]{Sh2}
H.\ Shimakura, The automorphism group of the vertex operator algebra $V_L^+$ for an even lattice $L$ without roots, {\it J. Algebra} {\bf 280} (2004), 29--57.

\end{thebibliography}
\end{document}